\documentclass[11pt]{amsart}
\usepackage{color}
\usepackage{url}

\renewenvironment{proof}[1][\proofname]{\par \normalfont \trivlist
\item[\hskip\labelsep\itshape #1]\ignorespaces
}{%
\hspace*{\fill}$\Box$ \endtrivlist }
\renewcommand{\proofname}{{\bf Proof}}
\newcommand{\romani}{\mathrm{i}}

\newcommand{\sgn}{\mbox{\textnormal{sgn\,}}}
\newcommand{\E}{\mbox{\textnormal{E}}}

\newcommand{\EP}{\E_\PP}
\newcommand{\PP}{\mathbb{P}}
\newcommand{\R}{\mathbb{R}}
\newcommand{\Q}{\mathbb{Q}}

\newcommand{\N}{\mathbb{N}}

\newcommand{\half}{\mbox{$\textstyle{\frac{1}{2}}$}}

\newtheorem{theorem}{Theorem}[section]
\newtheorem{lemma}[theorem]{Lemma}

\newtheorem{prop}[theorem]{Proposition}
\newtheorem{cor}[theorem]{Corollary}

\theoremstyle{definition}

\theoremstyle{remark}
\newtheorem{remark}[theorem]{Remark}

\numberwithin{equation}{section}

\begin{document}

\title{Negative volatility for a 2-dimensional square
root SDE}

\author{Peter Spreij}
\address{Korteweg-de Vries Instute of Mathematics, University of Amsterdam, Plantage Muidergracht 24, Amsterdam, The Netherlands}
\email{spreij@uva.nl}

\author{Enno Veerman}
\address{Korteweg-de Vries Instute of Mathematics, University of Amsterdam, Plantage Muidergracht 24, Amsterdam, The Netherlands}
\email{e.veerman@uva.nl}


\date{\today}



\begin{abstract}
In affine term structure models the short rate is modelled as an
affine transformation of a multi-dimensional square root process.
Sufficient conditions to avoid negative volatility factors are the
multivariate Feller conditions. {We will prove their necessity for
a 2-dimensional square root SDE in canonical form by presenting a
methodology based on measure transformations and {the trivial fact
that a random variable assumes negative values if it has negative
expectation. We exploit the property that solutions to square root
SDEs have expectations which solve a system of linear differential
equations}. As {an aside} we will present two proofs for the
martingale property of the density processes used in completely
affine models.}
\end{abstract}
\maketitle
\section{Introduction}
\subsection{Problem and motivation}

In recent years, affine term term structure models (ATSMs) have
become a popular instrument for modelling the dynamics of a term
structure, i.e.\ the dynamics of the short interest rate and the
long interest rate. These models have been introduced
by~\cite{dk96} and can be regarded as a multi-dimensional
extension of the Cox-Ingersoll-Ross model~\cite{cir85}. The short
rate is modelled as an affine transformation of a (possibly
multi-dimensional) state factor $X$ which satisfies a
multi-dimensional {square root SDE}. The diffusion part involves
square roots of affine transformations of $X$, which are called
{volatility factors}, or just {volatilities}. {Conditions need to
be imposed on the parameters to guarantee that the volatility
factors do not become negative, in order to assure pathwise
uniqueness and to justify the Feynman-Kac formula for the bond
price, see~\cite{eve06} for a detailed discussion. As shown
in~\cite{dk96} sufficient conditions for this are the so-called
\emph{multivariate Feller conditions}}, but they are not known to
be necessary (see~\cite[footnote 6]{ds00}).

Imposing the Feller conditions is not always desirable in
practice, as they might contradict with certain economic
principles. In~\cite{svv08} this is observed for a 2-factor model
where the state factor consists of the interest rate and
inflation. It turns out that estimating the model without the
Feller conditions {yields parameter values that are in agreement
with economic theory}. Therefore, the question is raised whether
these conditions can be relaxed. However, the model proposed
in~\cite{svv08} is actually a discrete-time ATSM. Imposing the
Feller conditions for excluding negative volatility factors is
meaningless in discrete time, since negative volatilities always
occur with positive probability, due to the normally distributed
jumps of the process. Instead, \cite{svv08} mainly investigates
the mathematical correctness of the discrete-time ATSM, with or
without the Feller conditions.

In the present paper we try to answer the initial question whether
the Feller conditions are necessary for excluding negative
volatility factors in continuous time. We focus on two
2-dimensional square root SDEs in canonical form, one with
proportional volatilities and {one} with {linear independent
volatility factors}, {see further down for a precise formulation}.
{The dimension is restricted to 2, since more or less only for
this case explicit computations can be performed.} The
proportional case is interesting from a practical point of view,
as it is the underlying SDE for one of the 2-factor models
proposed in~\cite{svv08}. Simulations suggest that for some
parameters, the volatility factor stays positive almost surely,
even though the Feller conditions do not hold. {We refer to
Figure~1 in~\cite{svv08} for an illustration.} However, in this
paper we prove that this suggestion is false: without the Feller
conditions there is always positive probability that the
volatility attains a negative value.
\subsection{Notation and definitions} With $\mathbb{S}_m(p)$ we denote the class of
$p$-dimensional square root SDEs with $m$ volatility factors. That
is, an element in $\mathbb{S}_m(p)$ is an SDE of the form
\begin{equation}\label{eq:smp}
dX_t=(aX_t+b)dt+\Sigma\sqrt{{|v(X_t)}|}dW_t,\quad X(0)=x_0\in\R^p,
\end{equation}
with {$X$ a $p$-dimensional stochastic process}, {$W$ a
$p$-dimensional Brownian motion, $a\in\R^{p\times p}$, $b\in\R^p$,
$\Sigma\in\R^{p\times p}$ non-singular. The $j$-th column of
$\Sigma$ is denoted by $\Sigma^j$. Furthermore, {$v(X_t)$ denotes
the diagonal matrix with diagonal elements
$V_{i,t}=v_i(X_t)=\alpha_i+\beta_i X_t$}, which we call
\emph{volatility factors} or \emph{volatilities}. Here
$\alpha_i\in\R$ and $\beta_i$ is a $p$-dimensional row-vector. We
let $\alpha=(\alpha_1,\ldots,\alpha_p)$, $\beta$ the matrix with
$i$-th row equal to $\beta_i$ and $m$ the rank of $\beta$}. {The
initial value $x_0$ is taken such that $v_i(x_0)\geq0$ for all
$i$.} {We write $\sqrt{|v(X_t)|}$ for the diagonal matrix with
diagonal elements $\sqrt{|V_{i,t}|}$}. When $V_{i,t}\geq0$ a.s.\
for all $i$, we omit the absolute sign and write $\sqrt{V_t}$
instead. {Later on we will use the notation $\sgn(v(X_t))$ for the
diagonal matrix with diagonal elements $\sgn(V_{i,t})$.}

{Existence of a weak solution to an SDE in $\mathbb{S}_m(p)$
holds, since the drift and diffusion part are continuous functions
which in addition fulfil a growth-condition, see Theorem~IV.2.3
and IV.2.4 in~\cite{Ikwat}.
To prove that existence and uniqueness of a \emph{strong} solution
holds (equivalent to \emph{pathwise uniqueness} by Theorem~IV.2.1
in~\cite{Ikwat}) appears very difficult for general square root
SDEs. The diffusion part is not Lipschitz-continuous, so standard
results, like Theorem~IV.3.1 in~\cite{Ikwat}, are not applicable.
Instead, one can use Theorem~1 in~\cite{YamadaI}, but this result
only applies to square root SDEs which can be written in a certain
canonical form (denoted by $\mathbb{A}_m(p)$ and $\mathbb{C}_p(p)$
below), for example when the Feller conditions are satisfied. In
general though, we do not know whether existence and uniqueness of
a strong solution holds for a square root SDE in
$\mathbb{S}_m(p)$.}

Sufficient conditions for strictly positive volatility factors
$V_{i,t}$ are the {multivariate Feller conditions
from~\cite{dk96}}, named after Feller's test for explosions. We
consider a weak version of these conditions which are sufficient
for non-negative instead of strictly positive volatility factors.
These \emph{weak Feller conditions} are given by
\begin{align}
\forall i,\forall j:&\,\beta_i\Sigma^j  = 0 \mbox { or $\partial\mathcal{D}_i\subset\partial\mathcal{D}_j$}, \label{eq:wfeller1} \\
\forall i,\forall x\in\partial\mathcal{D}_i:&\beta_i(a x+b) \geq0.
\label{eq:wfeller2}
\end{align}
Here $\partial\mathcal{D}_i$ denotes the boundary
$\{x\in\R^p:v_i(x)=0,v_j(x)\geq0,\forall j \}$.

The subclass ${\mathbb{C}_m(p)}\subset\mathbb{S}_m(p)$ contains
the square root SDEs which are in \emph{canonical form}, i.e.\
\begin{equation}\label{eq:can}
\Sigma=I,\quad V_i=X_i,\mbox{ for $i\leq m$,}
\end{equation}
with $I$ the identity matrix. We adopt the notation of~\cite{ds00}
for the class $\mathbb{A}_m(p)$, the SDEs which are in canonical
form and in addition \emph{satisfy the weak Feller conditions.} In
canonical form, the Feller conditions translate as
\begin{eqnarray}\label{eq:canfel}
\lefteqn{{i,j\leq m},k>m\Longrightarrow}  \\ & & a_{ij}\geq0
\mbox{ for $i\not=j$, } a_{ik}=0, b_i\geq 0,
\alpha_k\geq0,\beta_{ki}\geq0. \nonumber
\end{eqnarray}
If $X$ solves a square root SDE in $\mathbb{S}_m(p)$ which
satisfies the weak Feller conditions, then there exists an affine
transformation of $X$ which solves an SDE in canonical form, see
the appendix in~\cite{ds00} and Chapter 4 in~\cite{eve06}. It is
remarkable that when the volatilities are {proportional}, the
Feller conditions are not needed for this. This will be proved in
Proposition~\ref{prop:can}. We write
$\mathbb{S}(p){\subset\mathbb{S}_1(p)}$ for the $p$-dimensional
square root SDEs with proportional volatilities ({that is,
$\alpha_i=\alpha_1$ and $\beta_i=\beta_1$ for all $i$}), and
similarly $\mathbb{C}(p)$ for those in canonical form.
\subsection{Set-up}
{The rest of this paper is organized as follows. In
Section~\ref{sec:meth} we present the methodology to prove the
necessity of the Feller conditions for non-negative volatility
factors for the general class $\mathbb{S}_m(p)$. Explicit
computations are only possible for the special cases
$\mathbb{C}_2(2)$ and $\mathbb{S}(2)$, which we provide in the
remaining sections. The method is based on solving a system of
linear ODEs {satisfied by $\E X_t$} and transforming the
underlying probability measure via a certain exponential density
process $L$. For the method to work, it is necessary that $L$ is a
martingale. This is relatively straightforward for the class
$\mathbb{C}_2(2)$, but much more difficult to show for
$\mathbb{S}(2)$. Therefore, we first prove the necessity of the
Feller conditions for the class $\mathbb{C}_2(2)$ in
Section~\ref{sec:negvol2}, before tackling the harder case
$\mathbb{S}(2)$. We use two sections for working out the
methodology for the latter. Section~\ref{sec:nov} is entirely
devoted to verifying a local version of Novikov's condition (as
given in Corollary~3.5.14 in~\cite{Karshr}), in order to prove the
martingale property of $L$ for $\mathbb{S}(2)$.
Section~\ref{sec:negvol} deals with solving the systems of linear
ODEs and proving the necessity of the Feller conditions for
$\mathbb{S}(2)$.}

Though slightly off-topic, we have added an appendix with the
proof that $L$ is also {a} martingale for the class
$\mathbb{A}_m(p)$. We have two reasons for this. In
\emph{completely affine models} (see~\cite{duf02}) one uses this
particular exponential process $L$ to relate the physical with the
risk-neutral measure for SDEs in the class $\mathbb{A}_m(p)$.
However, the fact that $L$ is a legitimate density process (i.e.\
a martingale) is obscured in {the} literature. Therefore we
{clarify this once and for all}. We give two proofs that $L$ is a
martingale. Both serve as underlying {ideas} for proving the
martingale property of $L$ for the other classes $\mathbb{S}(2)$
and $\mathbb{C}_2(2)$, {which {is} the second reason.}
\section{Methodology}\label{sec:meth}
\noindent{This section presents the methodology to prove the
necessity of the (weak) Feller conditions with respect to
excluding negative volatility factors. The underlying idea applies
to the general class $\mathbb{S}_m(p)$, though explicit
computations are un-doable for higher dimensions.}
{The general scheme} for proving necessity of the Feller
conditions {consists of the following steps}:
\begin{enumerate}
\item {Let $(X,W)$ be a weak solution to~(\ref{eq:smp}) on some filtered
probability space $(\Omega,\mathcal{F},(\mathcal{F}_t),\PP)$.}
Then $\EP X_{t}$ solves a linear ODE by Lemma~\ref{lem:exode}
{below}, so we can compute $\EP V_{i,t}$.\label{enu:exode}
\item Let
\begin{equation}\label{eq:den}
L^\lambda_t:=\mathcal{E}(\int_0^\cdot\lambda^\top\sgn({v(X_s)})\sqrt{|{v(X_s)}|}d
W_s)_t,
\end{equation}
where $\lambda\in\R^p$. Fix an arbitrary time interval $[0,T]$,
{with $T>0$}. If the process $L^\lambda_t$ is a martingale on
$[0,T]$, then we can transform the measure $\PP$ into an
equivalent probability measure $\Q^\lambda$ on $\mathcal{F}_T$ by
$d\Q^\lambda=L^\lambda_Td\PP$.
By Girsanov's Theorem {(Theorem~3.51 in~\cite{Karshr})},
$W^\lambda$ defined by
$dW^\lambda_t=dW_t-\sgn({v(X_t)})\sqrt{|{v(X_t)}|}\lambda dt$ is a
Brownian motion under $\Q^\lambda$ on {$[0,T]$}. Moreover, the
resulting SDE for $X$ under $\Q^\lambda$ is still a square root
SDE:
\begin{equation}\label{eq:newSDE}
dX_t=(aX_t+b+\Sigma
v(X_t)\lambda){dt}+\Sigma\sqrt{|v(X_t)|}dW_t^\lambda,
\end{equation}
{where we can view the integral with respect to ${W}^\lambda$ as a
stochastic integral under $\Q^\lambda$ by Proposition~7.26
of~\cite{Jacod}.} {As in Step 1, also} $\E_{\Q^\lambda} X_{t}$
solves a linear ODE.\label{enu:mart}
\item Under violation of the Feller conditions,
for each $t>0$ we find $\lambda\in\R^p$ such that $\E_{\Q^\lambda}
V_{i,t}<0$. Then obviously $\Q^\lambda(V_{i,t}<0)>0$ and by
equivalence of measures also $\PP(V_{i,t}<0)>0$.\label{enu:neg}
\end{enumerate}
\begin{lemma}\label{lem:exode}
Let {$(X,W)$} be a weak solution {on some filtered probability
space $(\Omega,\mathcal{F},(\mathcal{F}_t),\PP)$} to
\[
dX_t=(a X_t+ b)dt+\sigma(t,X_t)dW_t,
\]
with $W$ a $p$-dimensional Brownian motion, $a\in \R^{p\times p}$,
$b\in\R^p$, $\sigma:[0,\infty)\times \R^p\rightarrow\R^{p\times
p}$ measurable and satisfying the {growth condition}
\[
\|\sigma(t,x)\|^2\leq K(1+\|x\|^2),\mbox{ for some positive
constant $K$}.
\]
If $\E\| X_0\|^2<\infty$, then $\bar{x}_t=\E X_t$ solves the ODE
\[
d\bar{x}_t=(a\bar{x}_t+b)dt,\quad \bar{x}(0)=\E X_0.
\]
\end{lemma}
\begin{proof}  Taking expectations gives
        \[
            \E X_t = \E X_0 + \E \int_0^t(a X_s+b)ds+\E
            \int_0^t\sigma(s,X_s)dW_s.
        \]
        By application of Problem
        5.3.15 in~\cite{Karshr} it holds that
        \[
            \E \max_{0\leq s\leq t}\|X_s\|^2 <\infty.
        \]
        In addition to the {growth condition} this implies that the stochastic integral is a martingale, whence its expectation equals
        zero. The result then follows by an application of Fubini.
\end{proof}
\medskip
{We apply this methodology to the classes $\mathbb{C}_2(2)$ and
$\mathbb{S}(2)$. For SDEs in $\mathbb{C}_2(2)$ pathwise uniqueness
holds, which we can use in {Step~\ref{enu:mart}} to verify the
martingale property of the exponential process $L^\lambda_t$.
However, for SDEs in {$\mathbb{S}(2)$} it is not clear whether we
have pathwise uniqueness. Therefore, we use an alternative method
to prove that $L^\lambda_t$ is a martingale by verifying Novikov's
condition. We {are able to} do this only up to the stopping time
$\tau=\inf\{t>0:V_{1,t}<0\}$ though. Consequently, we have to take
\begin{equation}\label{eq:den2}
L^\lambda_t:=\mathcal{E}(\int_0^\cdot\lambda^\top\sqrt{{v(X_{s\wedge\tau})}}d
W_s)_t,
\end{equation}
for the density process instead of (\ref{eq:den}). The result
$\PP(V_{i,t}<0)>0$ obtained in {Step~\ref{enu:neg}} needs to be
replaced by $\PP(\tau<t)>0$.}
\section{Negative volatility for $\mathbb{C}_2(2)$ without Feller
conditions}\label{sec:negvol2}\noindent{Consider the class
$\mathbb{C}_2(2)$ of square root SDEs with independent
volatilities ({meaning that} $\beta$ has full rank). For
notational convenience we write $V_t$ for the first coordinate of
a solution to the square root SDE, and $Y_t$ for the second
coordinate. So we consider SDEs of the form}
\begin{align}
d V_t &= (a_{11}V_t+a_{12} Y_t+b_1)dt + \sqrt{|V_t|}dW_{1,t},\qquad V_0=v_0\geq0,\label{eq:vv}\\
d Y_t &= (a_{21}V_t+a_{22} Y_t +b_2)dt +
\sqrt{|Y_t|}d{W}_{2,t},\qquad Y_0=y_0\geq0.\label{eq:yy}
\end{align}
{The Feller conditions~(\ref{eq:canfel}) in this case read
$a_{12},a_{21}\geq 0$ and $b_1,b_2\geq 0$. We shall violate the
condition $a_{12}\geq 0$ by assuming
\begin{equation}\label{eq:a12neg}
a_{12}<0,
\end{equation}
whereas we strengthen $a_{21}\geq 0$ and {$b_2\geq 0$} to}
\begin{equation}\label{eq:condition2}
a_{21}>0, {b_2>0}.
\end{equation}
{In this section we apply the method as described in
Section~\ref{sec:meth} to show that $V_t$ attains a negative value
with positive probability, for all $t>0$. The
SDE~(\ref{eq:newSDE}), obtained after the measure transformation
described in {Step~\ref{enu:mart}}, assumes the form
\begin{align}
d V_t &= (a_{11}^\lambda V_t+a_{12} Y_t+b_1)dt + \sqrt{|V_t|}dW_{1,t}^\lambda,\qquad V_0=v_0\geq0,\label{eq:vvnew}\\
d Y_t &= (a_{21}V_t+ a_{22}^\lambda Y_t +b_2)dt +
\sqrt{|Y_t|}d{W}_{2,t}^\lambda,\qquad
Y_0=y_0\geq0\label{eq:yynew},
\end{align}
{with $a_{11}^\lambda=a_{11}+\lambda_1$ and
$a_{22}^\lambda=a_{22}+\lambda_2$}. So in the corresponding ODEs
for the expectation, the parameters $a_{11}^\lambda$ and
$a_{22}^\lambda$ {depend on the chosen underlying probability
measure and are thus free to choose}. We show that these can be
{selected} in such a way that the first coordinate $\E V_t$ gets
negative, from which it follows that $V_t$ gets negative with
positive probability.} {Below we suppress the dependence of
$a_{11}$ and $a_{22}$ on $\lambda$}.
\begin{prop}\label{prop:solveODE2}
Let $a_{12}$, $a_{21}$, $b_1$, $b_2$, $x_0\geq 0$, $y_0\geq0$ be
arbitrary but fixed parameters and let $a_{11}$ and $a_{22}$ be
variable. Consider the family of systems of differential equations
parameterized by $a_{11},a_{22}$:
\begin{align}
\dot{x} & = a_{11}x+a_{12}y+b_1,\qquad x(0)=x_0\geq0;\\
\dot{y} & = a_{21}x+a_{22}y+b_2,\qquad y(0)=y_0\geq0.
\end{align}
Write $x(t,a_{11},a_{22})$ for the solution $x(t)$ depending on
$a_{11}$ and $a_{22}$. Assume~(\ref{eq:a12neg})
and~(\ref{eq:condition2}). Then for all $t_0>0$ there exist
$a_{11}$ and $a_{22}$ such that  $x(t_0,a_{11},a_{22})<0$.
\end{prop}
\begin{proof} {We use the following notation: $\tau$ is the trace of $a$,
$\Delta$ {its} determinant, $\rho=a_{12}b_2-a_{22}b_1$,
$\bar{x}=\rho/\Delta$, {$D=\tau^2-4\Delta$}. By eliminating $y$}
we obtain a second order equation for $x$:
\begin{equation}
\ddot{x}-\tau \dot{x} +\Delta x -\rho=0.
\end{equation}
If $\Delta\not=0$ this has the general solution
\[
x(t)=B_1e^{r_1 t}+B_2 e^{r_2 t}+\bar{x},
\]
where $r_i=\half(\tau\pm\sqrt{D})$, {$i=1,2$}. We take $a_{11}=0$
and $a_{22}>0$ such that $D>0$. Notice that {$a_{11}=0$ implies}
$\Delta>0$. {Solving {for} $B_2$ gives
\begin{align*}
B_2&=\frac{r_1(\bar{x}-x_0)+a_{12}y_0+b_1}{r_2-r_1}=\frac{a_{12}y_0}{a_{22}}+\frac{a_{12}(b_2+a_{21}x_0)}{a_{22}^2}+O(a_{22}^{-3}),
\end{align*}
as $a_{22}\rightarrow\infty$. Hence $B_2<0$ for $a_{22}$ big
enough by the assumptions~(\ref{eq:a12neg})
and~(\ref{eq:condition2}). Furthermore it holds that
\[
r_1=O(a_{22}^{-1}),\,r_2=a_{22}-O(a_{22}^{-1}),\mbox{ as
$a_{22}\rightarrow\infty$,}
\]
From this it easily follows that for arbitrary $t_0>0$ we have
\begin{align*}
{x(t_0,0,a_{22})}&=B_1e^{r_1 t_0}+B_2 e^{r_2
t_0}+\bar{x}\\&=x_0e^{r_1 t_0}-\bar{x}(e^{r_1 t_0}-1)+B_2(e^{r_2
t_0}-e^{r_1 t_0})\rightarrow-\infty,
\end{align*}
as $a_{22}\rightarrow\infty$, since $B_2e^{r_2 t_0}$ tends to
$-\infty$ and dominates the other terms.} {Hence, the choice
$a_{11}=0$ and $a_{22}>0$ big enough results in
$x(t_0,a_{11},a_{22})<0$}.
\end{proof}
{
\begin{theorem}\label{theorem:PP}
{Let $((V,Y),W)$ be a solution to (\ref{eq:vv}), (\ref{eq:yy}), on
some filtered probability space
$(\Omega,\mathcal{F},(\mathcal{F}_t),\PP)$}.
{Assume~(\ref{eq:condition2}) and that the Feller conditions are
violated by~(\ref{eq:a12neg})}. Then for all $t>0$ it holds that
\[
\PP(V_t<0)>0.
\]
\end{theorem}
\begin{proof} We follow the methodology as described in Section~\ref{sec:meth}. Time is restricted to an arbitrary but finite
interval $[0,T]$, {with $T>0$}. In addition to the existence of a
weak solution we have pathwise uniqueness by Theorem~1
in~\cite{YamadaI}, {which implies existence of a strong solution
by Theorem~IV.2.1 in~\cite{Ikwat}.} Hence we can apply
Proposition~\ref{prop:change2} to obtain that {$L_t^\lambda$ as
defined by~(\ref{eq:den})} is a martingale for all
$\lambda\in\R^2$. So we can change the measure $\PP$ by
$d\Q^\lambda=L_T^\lambda d\PP$ and obtain an SDE under
$\Q^\lambda$ {up to time $T$}, {as given by~(\ref{eq:vvnew})
and~(\ref{eq:yynew})}{.}
By Lemma~\ref{lem:exode} and Proposition~\ref{prop:solveODE2}, for
all {$t\in(0,T]$} we can choose $\lambda$ such that
$\E_{\Q^\lambda} V_{t}<0$, which implies that
$\Q^\lambda(V_t<0)>0$. By equivalence of measures it follows that
$\PP(V_t<0)>0$ for all {$t\in(0,T]$}. {Since $T>0$ was chosen
arbitrarily, the result holds for all $t>0$.}
\end{proof}
}
\section{Measure transformation {for} {$\mathbb{S}(2)$} without Feller conditions}\label{sec:nov}
\noindent{Now that we have applied the methodology of
Section~\ref{sec:meth} to the class $\mathbb{C}_2(2)$, we {would
like} to do the same for $\mathbb{S}(2)$. Since we do not know
whether pathwise uniqueness holds for {solutions to SDEs from}
this class, we need to do some extra work in {Step~\ref{enu:mart}}
for verifying the martingale property of $L_t^\lambda$. This is
done in the current section by checking a local version of
Novikov's condition. Then we work out the remaining steps for
proving necessity of the Feller conditions {in the next section}.
First, however, we show that \emph{every} SDE in $\mathbb{S}(p)$
can be rewritten {in canonical form}.
\begin{prop}\label{prop:can} {Let $(X,W)$ be a
solution on some filtered probability space
$(\Omega,\mathcal{F},(\mathcal{F}_t),\PP)$} of a $p$-dimensional
square root SDE from $\mathbb{S}(p)$ with one volatility factor
$V_{1,t}=\alpha_1+\beta_1 X_t$ ($\alpha_1\in\R$, $\beta_1$ a
$p$-dimensional row vector, not equal to zero):
\[
dX_t = (aX_t+b)dt+\Sigma\sqrt{|V_{1,t}|}dW_t,
\]
$a\in\R^{p\times p}$, $b\in\R^{p}$ and $\Sigma\in\R^{p\times p}$
non-singular. Then there exists an affine transformation
$\widetilde{X}$ of $X$ such that $\widetilde{X}$ solves an SDE
from $\mathbb{C}(p)$:
\begin{equation}\label{eq:ortho}
d\widetilde{X}_t =
(\widetilde{a}\widetilde{X}_t+\widetilde{b})dt+\sqrt{|\widetilde{X}_{1,t}|}d\widetilde{W}_t,
\end{equation}
where $\widetilde{W}$ is an orthogonal transformation of $W$,
whence also a Brownian motion. Moreover,
$\widetilde{X}_{1,t}=cV_{1,t}$ for some positive constant $c>0$.
In addition we can take $\widetilde{a}_{1j}\geq0$ for $j\not=1$.
\end{prop}}
\begin{proof} We need to find $K\in\R^{p\times p}$ and
$\ell\in\R^p$ such that for $\widetilde{X}=KX+\ell$ we have
\[
\sqrt{|V_{1,t}|}K\Sigma
dW_t=\sqrt{|\widetilde{X}_{1,t}|}d\widetilde{W}_t,
\]
i.e.\ $\widetilde{X}_{1,t}=cV_{1,t}=c\alpha_1+c\beta_1 X_t$ for
some $c>0$ and $K\Sigma/\sqrt{c}$ is orthonormal. The first
requirement is fulfilled if $K_1=c\beta_1$ and $\ell_1=c\alpha_1$.
For the second requirement we need that the first row of
$K\Sigma/\sqrt{c}$ has length one, so
$\|(K\Sigma/\sqrt{c})_1\|=\|K_1\Sigma/\sqrt{c}\|=\|\sqrt{c}\beta_1\Sigma\|=1$,
i.e.\ $c=1/\|\beta_1\Sigma\|^2$. This gives that
$\ell_1=c\alpha_1=\alpha_1/\|\beta_1\Sigma\|^2$ and we may take
$\ell_i$ arbitrarily for $i\not=1$. Moreover the remaining row
vectors of $K\Sigma/\sqrt{c}$ should be mutually orthogonal and
orthogonal to the first row vector, while also be of length one.
Choose such vectors $k_j$, $j=2,\ldots,p$ and write $M$ for the
matrix with these row vectors. Then we take
\[
K=\sqrt{c}
  \begin{pmatrix}
    \beta_1\Sigma/\|\beta_1\Sigma\| \\
    M
  \end{pmatrix}\Sigma^{-1}
    =\begin{pmatrix}
    \beta_1/\|\beta_1\Sigma\|^2 \\
    M\Sigma^{-1}/\|\beta_1\Sigma\|
  \end{pmatrix}.
\]
Thus the SDE for $\widetilde{X}=KX+\ell$ is of the canonical
form~(\ref{eq:ortho}). If $\widetilde{a}_{1j}<0$ for $j\not=1$,
then we take $-\widetilde{X}_{j,t}$ instead of
$\widetilde{X}_{j,t}$, which also gives a canonical SDE (replacing
$\widetilde{W}_{j,t}$ by $-\widetilde{W}_{j,t}$ still gives a
Brownian motion), but with $\widetilde{a}_{1j}>0$ instead.
\end{proof}
\medskip
{As a consequence of the previous proposition, it is sufficient to
consider necessity of the Feller conditions~(\ref{eq:canfel}) for
the SDEs in {$\mathbb{C}(2)$}, {which read $a_{12}=0$ and
$b_1\geq0$. We only consider the first condition.} In view of
Proposition~\ref{prop:can} we may always assume $a_{12}\geq 0$,
{so we violate the first Feller condition by}
\begin{equation}\label{eq:a12pos}
a_{12}>0.
\end{equation}
}
\begin{remark} When $a_{12}=0$ we are essentially dealing with a
1-dimensional square root SDE. Then $b_1\geq0$ is the remaining
Feller condition, for which necessity follows by the 1-dimensional
Feller's test for explosions (Theorem~5.5.29 in~\cite{Karshr}).
\end{remark}
\noindent{As in Section~\ref{sec:negvol2}}, we write $V_t$ for the
volatility factor $X_{1,t}$. {Moreover, for $a_{12}> 0$ we may}
substitute $Y_t=X_{2,t}+b_1/a_{12}$ for the second coordinate so
that the resulting SDE is of the form
\begin{align}
d V_t &= (a_{11}V_t+a_{12} Y_t)dt + \sqrt{|V_t|}dW_{1,t},\qquad\qquad V_0=v_0\geq0,\label{eq:v}\\
d Y_t &= (a_{21}V_t+a_{22} Y_t +b_2)dt +
\sqrt{|V_t|}d{W}_{2,t},\qquad Y_0=y_0\in\R.\label{eq:y}
\end{align}
{Hence, for $a_{12}> 0$ we can assume without loss of generality
that $b_1=0$. Note that in the present notation we have
$\tau=\inf\{t>0:V_t<0\}$.}

We prove that $L_t^\lambda$ defined by~(\ref{eq:den2}) is a
martingale for all $\lambda\in\R^2$ by verifying Novikov's
condition. The proof is similar to {the proof of}
Proposition~\ref{prop:nov}, but more complicated due to the
violation of the Feller conditions.
{Notice} that for $\mathbb{A}_m(p)$, Novikov's condition is
satisfied under the additional requirement~(\ref{eq:addreq}).
Likewise, for $\mathbb{C}(2)$ without Feller conditions we need an
additional condition to justify Novikov's condition, namely
\begin{equation}\label{eq:condition}
\mbox{$a_{11}<0$, $a_{22}<0$ and $\det a>0$.}
\end{equation}
Note that this implies negative real parts for the eigenvalues of
the matrix $a$, whence a solution $X$ is mean-reverting. This
latter property would suggest more ``stability'' for $X$, and
hence more integrability properties (like Novikov's condition).

We first prove some lemmas.
\begin{lemma}\label{lemV}
{Let $((V,Y),W)$ be a solution to~(\ref{eq:v}), (\ref{eq:y}), on
some filtered probability space
$(\Omega,\mathcal{F},(\mathcal{F}_t),\PP)$}. Assume $a_{11}<0$.
For $0\leq c\leq -a_{11}$ it holds that
\[
\E \exp(cV_{t\wedge\tau})\leq\exp(cv_0)
\left[\E\exp\left(2ca_{12}\int_0^{t\wedge\tau}Y_sds\right)\right]^{1/2}.
\]
\end{lemma}
\begin{proof} From~(\ref{eq:v}) {one obtains}
\begin{align*}
\E \exp(cV_{t\wedge\tau})
&=\exp(cv_0)\E\Big[\exp\left(c\int_0^{t\wedge\tau}\sqrt{|V_s|}dW_{1,s}+ca_{11}\int_0^{t\wedge\tau}V_sds\right)\\
&\qquad\qquad\qquad\times\exp\left(ca_{12}\int_0^{t\wedge\tau}Y_sds\right)\Big]\\
&\leq\exp(cv_0)\left[\E\exp\left(2c\int_0^{t\wedge\tau}\sqrt{|V_s|}dW_{1,s}+2ca_{11}\int_0^{t\wedge\tau}V_sds\right)\right]^{1/2}\\
&\qquad\qquad\qquad\times\left[\E\exp\left(2ca_{12}\int_0^{t\wedge\tau}Y_sds\right)\right]^{1/2},
\end{align*}
where the last inequality follows from Cauchy-Schwarz. For $0\leq
c\leq -a_{11}$ it holds that
\begin{align*}
&\exp\left(2c\int_0^{t\wedge\tau}\sqrt{|V_s|}dW_{1,s}+2ca_{11}\int_0^{t\wedge\tau}V_sds\right)
\leq\mathcal{E}(\int_0^\cdot
2c\sqrt{|V_s|}dW_{1,s})_{t\wedge\tau},
\end{align*}
since $2c(c+a_{11})\leq 0$ and $V_s\geq0$ holds for all
$s\leq\tau$. By optional sampling (see for example Problem 1.3.23
in~\cite{Karshr}), the stopped exponential process in the above
display is also a supermartingale. So it has expectation less than
or equal to 1 and the result follows.
\end{proof}
{
\begin{lemma}\label{lemZ}
Let $f:[0,\infty)\rightarrow\R$ be integrable. For all $t\geq0$,
$\varepsilon>0$ it holds that
\begin{align*}
\exp(\int_t^{t+\varepsilon} f(s) ds)&\leq
\frac{1}{\varepsilon}\int_t^{t+\varepsilon}\exp(\varepsilon
f(s))ds.
\end{align*}
\end{lemma}
\begin{proof} Fix $t\geq0$, $\varepsilon>0$. Define a probability
measure $\mu$ on the Borel sigma-algebra $\mathcal{B}(\R)$ by
\[
d\mu=\frac{1}{\varepsilon}1_{[t,t+\varepsilon]}d\lambda,
\]
with $\lambda$ the Lebesgue-measure. The exponential function is
convex, so we can apply Jensen's inequality to obtain
\begin{align*}\exp(\int_t^{t+\varepsilon}
f(s)ds)&=\exp\left(\int \varepsilon f(s)\mu(ds)\right) \leq \int
\exp(\varepsilon
f(s))\mu(ds)\\&=\varepsilon^{-1}\int_t^{t+\varepsilon}
\exp(\varepsilon f(s))ds.
\end{align*}
\end{proof}}
\begin{lemma}\label{lemY}
{Let $((V,Y),W)$ be a solution to~(\ref{eq:v}), (\ref{eq:y}), on
some filtered probability space
$(\Omega,\mathcal{F},(\mathcal{F}_t),\PP)$}. Assume $a_{22}< 0$,
$\det a>0$. Write
\begin{align*}
c_1=\frac{-2a_{22}\det a }{a_{12}^2+a_{22}^2},\quad
c_2=\frac{2a_{12}\det a}{a_{12}^2+a_{22}^2}.
\end{align*}
Then for $0\leq c\leq c_2$ it holds that
\[
\E \exp(c Y_t1_{{t\leq\tau}})\leq 1+\exp(k(t)),
\]
{with $k(t)=c_1 v_0+c_2 y_0+c_2b_2t1_{\{b_2>0\}}$.}
\end{lemma}
\begin{proof} Since $c_1a_{12}+c_2a_{22}=0$ and
$c_1a_{11}+c_2a_{21}=-\half(c_1^2+c_2^2)$, (\ref{eq:v})
and~(\ref{eq:y}) {give}
\begin{align*}
\E \exp(c_1&V_{t\wedge\tau}+c_2Y_{t\wedge\tau})=\\&=\E
\exp\Big(c_1v_0+c_2y_0+c_2b_2(t\wedge\tau)+\int_0^{t\wedge\tau}(c_1a_{11}+c_2a_{21})V_s ds\\
&\qquad\quad+\int_0^{t\wedge\tau}(c_1a_{12}+c_2a_{22})Y_s ds+\int_0^{t\wedge\tau}\sqrt{|V_s|}(c_1\,c_2) dW_s\Big)\\
&\leq
\exp({k(t)})\E\exp\Big(\int_0^{t\wedge\tau}\sqrt{|V_s|}(c_1\,c_2)  dW_s-\half\int_0^{t\wedge\tau}(c_1^2+c_2^2)V_s ds\Big)\\
&=
\exp({k(t)})\E\mathcal{E}(\int_0^\cdot\sqrt{|V_s|}(c_1\,c_2)dW_s)_{t\wedge\tau}\\
&\leq \exp({k(t)}).
\end{align*}
Note that $c_1\geq0$ and $c_2\geq0$, so for $0\leq c\leq c_2$ we
have
\begin{align*}
\E
\exp(cY_t1_{\{t\leq\tau\}})&\leq\E\exp(cY_t1_{\{t\leq\tau\}}1_{\{Y_t>0\}})\leq\E\exp(c_2Y_t1_{\{t\leq\tau\}}1_{\{Y_t>0\}})\\&\leq1+
\E\exp(c_2Y_{t\wedge\tau})\\
&\leq1+\E\exp(c_1V_{t\wedge\tau}+c_2Y_{t\wedge\tau}),
\end{align*}
which gives the result. \end{proof}
\begin{prop}\label{prop}
Let $((V,Y),W)$ be a solution to~(\ref{eq:v}), (\ref{eq:y}), on
some filtered probability space
$(\Omega,\mathcal{F},(\mathcal{F}_t),\PP)$.
{Assume~(\ref{eq:a12pos}) and (\ref{eq:condition}).} Fix an
arbitrary $c>0$ and define
\[
\varepsilon(t)=\min\left(-\frac{a_{11}}{c},\half\left(-t+\sqrt{t^2+\frac{2c_2}{ca_{12}}}\right)\right),
\]
with $c_2$ as in Lemma~\ref{lemY}. Then for all $t\geq0$ it holds
that
\[
\E\exp\left(c\int_t^{t+\varepsilon(t)}V_{s\wedge\tau}ds\right)<\infty.
\]

\end{prop}
\begin{proof} Fix $t\geq0$ and $\varepsilon:=\varepsilon(t)$. Applying
respectively Lemmas~\ref{lemZ},~\ref{lemV},~\ref{lemZ}
and~\ref{lemY}, we obtain
\begin{align*}
&\,\quad\E\exp\left(c\int_t^{t+\varepsilon}V_{s\wedge\tau}ds\right)\\&\stackrel{\ref{lemZ}}{\leq}
\frac{1}{\varepsilon}\int_t^{t+\varepsilon}\E\exp(\varepsilon c
V_{s\wedge\tau})ds\\
&\stackrel{\ref{lemV}}{\leq} \frac{1}{\varepsilon}\exp(\varepsilon
cv_0) \int_t^{t+\varepsilon}ds\left[\E\exp\left(2\varepsilon
ca_{12}\int_0^sY_u1_{\{u\leq\tau\}}du\right)\right]^{1/2}\\
&\stackrel{\ref{lemZ}}{\leq}\frac{1}{\varepsilon}\exp(\varepsilon
cv_0)\int_t^{t+\varepsilon} ds \left[\frac{1}{s}\int_0^{s}\E\exp(s
2\varepsilon
c a_{12} Y_u1_{\{u\leq\tau\}})du\right]^{1/2}\\
&\stackrel{\ref{lemY}}{\leq}\frac{1}{\varepsilon}\exp(\varepsilon
cv_0)\int_t^{t+\varepsilon} ds
\left[\frac{1}{s}\int_0^{s}(1+\exp({k(u)}))du\right]^{1/2}\\
&<\infty,
\end{align*}
{with $k(u)=c_1 v_0+c_2 y_0+c_2b_2u1_{\{b_2>0\}}$}. Note that to
apply Lemma~\ref{lemV} it is necessary that
\[
0\leq\varepsilon c\leq -a_{11},
\]
which holds true by definition of $\varepsilon(t)$. To apply
Lemma~\ref{lemY} we need to check that
\[
0\leq s2\varepsilon c a_{12}\leq c_2, \textnormal{ for all }t\leq
s\leq t+\varepsilon.
\]
Choosing $s=t+\varepsilon$ this comes down to
\[
\varepsilon^2+t\varepsilon - \frac{c_2}{2c a_{12}}\leq0.
\]
This is satisfied if and only if
\[
\half\left(-t-\sqrt{t^2+\frac{2c_2}{ca_{12}}}\right)\leq\varepsilon\leq\half\left(-t+\sqrt{t^2+\frac{2c_2}{ca_{12}}}\right),
\]
which indeed holds true. \end{proof}
\medskip
Using this proposition we can now verify the local version of
Novikov's condition, {as given in} Corollary~3.5.14
in~\cite{Karshr}.
\begin{prop}\label{prop:locnov}
Let $((V,Y),W)$ be a solution to (\ref{eq:v}), (\ref{eq:y}), on
some filtered probability space
$(\Omega,\mathcal{F},(\mathcal{F}_t),\PP)$.
{Assume~(\ref{eq:a12pos}) and (\ref{eq:condition}).} Then for all
$c>0$ there exist $0=t_0<t_1<t_2<\ldots<t_n\uparrow\infty$ such
that
\[
\E \exp( c\int_{t_i}^{t_{i+1}} V_{s\wedge\tau}
ds)<\infty,\textnormal{ for all }i,
\]
and
$L_t^\lambda=\mathcal{E}(\int_0^\cdot\lambda^\top\sqrt{V_{s\wedge\tau}}dW_s)_t$
is a martingale for all $\lambda\in\R^2$.
\end{prop}
\begin{proof} Fix $c>0$ and take $t_0=0$, $t_{i+1}=t_i+\varepsilon(t_i)$,
with $\varepsilon$ defined as in Proposition~\ref{prop}. The
result follows upon noting that
$t_n=\sum_{i=0}^{n-1}\varepsilon(t_i)\uparrow\infty$. The latter
can be proved by contradiction. Suppose $t_n\uparrow M<\infty$.
Then on the one hand $\varepsilon(t_n)\rightarrow 0$, since the
sum $\sum_{i=0}^{n-1}\varepsilon(t_i)$ converges. But on the other
hand,
\begin{align*}
\varepsilon(t_n)&=\min\left(-\frac{a_{11}}{c},\half\left(-t_n+\sqrt{t_n^2+\frac{2c_2}{ca_{12}}}\right)\right)\\
&\rightarrow
\min\left(-\frac{a_{11}}{c},\half\left(-M+\sqrt{M^2+\frac{2c_2}{ca_{12}}}\right)\right)>0,
\end{align*}
which is a contradiction. Hence $t_n\uparrow\infty$. Since the
local version of Novikov's condition holds true, $L_t^\lambda$ is
a martingale for all $\lambda\in\R^2$.
\end{proof}
\section{Negative volatility for $\mathbb{S}(2)$ without Feller conditions}\label{sec:negvol}
\noindent{In the previous section we showed that the measure
transformation given in Step~\ref{enu:mart} is legitimate for
$\mathbb{S}(2)$. Now we complete the proof of {the} necessity of
the Feller conditions by completing Step~\ref{enu:neg}. {Recall
that by Proposition~\ref{prop:can} we can write an SDE from
$\mathbb{S}(2)$ in canonical form {as given by~(\ref{eq:v})
and~(\ref{eq:y})}. Then the} SDE~(\ref{eq:newSDE}), obtained after
the measure transformation, assumes the form
\begin{align}
d V_t &= ({a_{11}^\lambda}V_t+a_{12} Y_t)dt + \sqrt{|V_t|}dW_{1,t}^\lambda,\qquad\qquad V_0=v_0\geq0,\label{eq:vnew}\\
d Y_t &= ({a_{21}^\lambda}V_t+a_{22} Y_t +b_2)dt +
\sqrt{|{V_t}|}d{W}_{2,t}^\lambda,\qquad
Y_0=y_0\in\R,\label{eq:ynew}
\end{align}
with $t<\tau$ and {$a_{11}^\lambda=a_{11}+\lambda_1$ and
$a_{21}^\lambda=a_{21}+\lambda_2$}. So in the corresponding ODEs
for the expectation, the parameters $a_{11}$ and $a_{21}$ {depend
on the chosen underlying probability measure and are thus free to
choose.} Analogously to Proposition~\ref{prop:solveODE2}, we
prove:}
\begin{prop}\label{prop:solveODE}
{{Let $a_{12}$, $a_{22}$, $b_2$, $x_0\geq 0$, $y_0$ be arbitrary
but fixed parameters  and let $a_{11}$ and $a_{21}$ be variable.}
Consider the family of systems of differential equations
{parameterized by $a_{11},a_{21}$}:}
\begin{align}
\dot{x} & = a_{11}x+a_{12}y,\qquad\qquad x(0)=x_0\geq0;\\
\dot{y} & = a_{21}x+a_{22}y+b_2,\qquad y(0)=y_0.
\end{align}
Write $x(t,a_{11},a_{21})$ for the solution $x(t)$ depending on
$a_{11}$ and $a_{21}$. If $a_{12}\not=0$ and
$(x_0,\dot{x}_0,\dot{y}_0)\not=(0,0,0)$, then it holds that for
all $t_0>0$ there exist $a_{11}$ and $a_{21}$ such that
$x(t_0,a_{11},a_{21})<0$.
\end{prop}
\begin{proof}
We use the same notation as in the proof of
Proposition~\ref{prop:solveODE2}, {but for reasons of brevity we
write $x(t)$ instead of $x(t,a_{11},a_{21})$}. Again, by
eliminating $y$ we obtain a second order equation for $x$:
\begin{equation}
\ddot{x}-\tau \dot{x} +\Delta x -\rho=0,
\end{equation}
If $D=\tau^{2}-4\Delta <0$, then the characteristic equation $
r^{2}-\tau r +\Delta =0 $ has two different complex roots, which
are $r_i=\half(\tau\pm \romani\sqrt{|D|})$. In that case the
differential equation for $x$ has the general solution
\[
x(t)=\exp(\half\tau t)(c_{1}\cos (\omega t) + c_{2}\sin (\omega
t))+\bar{x},
\]
with {$\omega=\half\sqrt{|D|}$} and $c_{1}$, $c_{2}$ are
determined by the initial conditions $x_{0}$ and $y_{0}$ of the
original system. Solving for $c_{1}$ and $c_{2}$ yields
\begin{align*}
c_{1} & = x_{0}-\bar{x}, \\
c_{2} & = \frac{1}{\omega}(\dot{x}_0 -\half\tau(x_{0}-\bar{x})).
\end{align*}
Without loss of generality we may assume $a_{12}>0$ as we can
substitute $-y$ for $y$ to change the sign of $a_{12}$. Note that
\begin{equation}\label{eq:atoinf}
a_{21}\to
-\infty\Longrightarrow\Delta\rightarrow\infty\Longrightarrow
D\rightarrow-\infty, \bar{x}\rightarrow 0\mbox{ and }
\omega\rightarrow\infty.
\end{equation}
Fix $t_0>0$ and suppose $x_0>0$. By~(\ref{eq:atoinf}) we can
choose $a_{21}$ such that $D<0$, $\omega =(\pi+2\pi k)/t_0$ for
some $k\in\N$, and $\bar{x}<(x_0 \exp(\half \tau t_0)/(\exp(\half
\tau t_0)+1)$. It follows that
\[
x(t_0)=-\exp(\half\tau t_0)c_{1}+\bar{x}=-x_0\exp(\half\tau
t_0)+(\exp(\half\tau t_0)+1)\bar{x}<0.
\]
If $x_0=0$ and $\rho\not=0$ then we take $a_{11}$ such that
$\sgn(\tau)=\sgn(\rho)$ and $a_{21}$ such that $D<0$, $\omega
=2\pi k/t_0$ for some $k\in\N$, which is possible in view
of~(\ref{eq:atoinf}). Then $\Delta>0$ and
$\sgn(\bar{x})=\sgn(\tau)$, so
\[
x(t_0)=\exp(\half\tau t_0)c_{1}+\bar{x}=(1-\exp(\half\tau
t_0))\bar{x}<0.
\]
If $x_0=\rho=0$ and $\dot{x}_0\not=0$, then we choose $a_{21}$
such that $D<0$, $\omega=(\half\pi+\pi k)/t_0$ with $k\in2\N$ if
$\dot{x}_0<0$ and $k\in2\N+1$ if $\dot{x}_0>0$. Then
\[
x(t_0)=\exp(\half\tau t)c_{2}\sin (\omega t_0))=\exp(\half\tau
t)\frac{\dot{x}_0}{\omega}\cdot(1_{\{\dot{x}_0<0\}}-1_{\{\dot{x}_0>0\}})<0.
\]
If $x_0=\rho=\dot{x}_0=0$ then $\dot{y}_0=0$, so this case is
excluded by assumption.
\end{proof}
\medskip
{Since Novikov's condition is only verified
under~(\ref{eq:condition}), we do not know whether $L_t^\lambda$
defined by~(\ref{eq:den2}) is a martingale without this condition.
Therefore, to apply the methodology of Section~\ref{sec:meth} to
the more general case, we need to do some extra work. First we
show in Theorem~\ref{theorem:P} the necessity of the Feller
condition $a_{12}=0$, when~(\ref{eq:condition}) does hold. Then in
Theorem~\ref{theorem:Q} we relax~(\ref{eq:condition}) to
$a_{22}<0$. We show that if $a_{22}<0$ and the Feller condition
$a_{12}=0$ is violated, there exists a solution $((V,Y),W)$ to the
SDE~(\ref{eq:vnew}), (\ref{eq:ynew}), up to the stopping time
$\tau$, such that $V$ gets negative with positive probability.
This is done as follows.}

We construct a solution to the SDE by first changing the other
parameters $a_{11}$ and $a_{21}$ in such a way
that~(\ref{eq:condition}) does hold, and obtaining a solution to
the corresponding SDE under some measure {$\Q$} (for which we know
that $V$ gets negative with positive {$\Q$}-probability from
Theorem~\ref{theorem:P}). Then changing the measure {$\Q$} into an
equivalent measure {$\PP$}, we retrieve the original SDE using
Girsanov's Theorem. By equivalence of measures $V$ will also get
negative under {$\PP$}.
\begin{theorem}\label{theorem:P}
{Let $((V,Y),W)$ be a solution to (\ref{eq:v}), (\ref{eq:y}), on
some filtered probability space
$(\Omega,\mathcal{F},(\mathcal{F}_t),\PP)$}. Assume the Feller
conditions are violated by~(\ref{eq:a12pos}). In addition
assume~(\ref{eq:condition}). Then for all $T>0$ it holds that
\[
\PP(\tau<T)>0.
\]
\end{theorem}
\begin{proof}
{We follow the methodology as described in Section~\ref{sec:meth}
and give a proof by contradiction. Time is restricted to an
arbitrary but finite interval $[0,T]$, {with $T>0$}.
Proposition~\ref{prop:locnov} gives that $L_t^\lambda$, as defined
by~(\ref{eq:den2}), is a martingale for all $\lambda\in\R^2$. So
we can change the measure on $\mathcal{F}_T$ by $d\Q^\lambda =
L^\lambda_T d\PP$ and obtain an SDE under $\Q^\lambda$, as given
by~(\ref{eq:vnew}) and~(\ref{eq:ynew}), {up to $\tau\wedge T$}.
Now assume $\PP(\tau<T)=0$. Then we can apply
Lemma~\ref{lem:exode}. By Proposition~\ref{prop:solveODE}, for all
{$t\in(0,T)$} we can choose $\lambda$ such that
$\E_{\Q^\lambda}V_t < 0$, which implies that
$\Q^\lambda(V_{t}<0)>0$ and by equivalence of measures also
$\PP(V_{t}<0)>0$, which contradicts the assumption that
$\PP(\tau<T)=0$.}
%
\end{proof}
{
\begin{theorem}\label{theorem:Q}
Consider an SDE in $\mathbb{C}(2)$ given by~(\ref{eq:v})
and~(\ref{eq:y}). Assume the Feller conditions are violated by
$a_{12}>0$. In addition assume $a_{22}<0$. Let time be restricted
to an arbitrary but finite interval $[0,T]$, {with $T>0$}. Then
there exists {an adapted stochastic process} $((V,Y),W)$ on some
filtered probability space
{$(\Omega,\mathcal{F},(\mathcal{F}_t)_{0\leq t\leq T},\PP)$ which
is a solution to~(\ref{eq:v}),~(\ref{eq:y}),} up to {$\tau\wedge
T$}, such that $\PP(\tau<t)>0$ for all {$t\in(0,T]$}.
\end{theorem}
\begin{proof}
Take $\lambda\in\R^2$ such that~(\ref{eq:condition}) holds true
with $a_{i1}$ replaced by $a_{i1}+\lambda_i$, for $i=1,2$. Let
$((V,Y),W^\lambda)$ be a weak solution
to~(\ref{eq:vnew}),~(\ref{eq:ynew}), on some filtered probability
space $(\Omega,\mathcal{F},(\mathcal{F}_t),\Q^\lambda)$, {with
time unrestricted}. For this the conditions of
Theorem~\ref{theorem:P} and Proposition~\ref{prop:locnov} hold
true. Applying the theorem gives that $\Q^\lambda(\tau<t)>0$ for
all $t$. Applying the proposition with $-\lambda$ instead of
$\lambda$ gives that
\[
L^{-\lambda}_t:=\mathcal{E}(\int_0^\cdot-\lambda^\top\sqrt{V_{s\wedge\tau}}d
W^\lambda_s)_t
\]
is a martingale. It follows that $\PP$ defined by
$d\PP=L_T^{-\lambda}d\Q^\lambda$ is a probability measure on
$\mathcal{F}_T$ equivalent to $\Q^\lambda$. Moreover, the process
$W$ defined by $dW_t=dW^\lambda_t+\sqrt{V_{t\wedge\tau}}\lambda
dt$ is a Brownian Motion on $[0,T]$ under $\PP$ by Girsanov's
Theorem and $((V,Y),W)$ solves
\begin{align*}
d V_t &= ({a_{11}^\lambda}V_t+a_{12}
Y_t-\lambda_1V_{t\wedge\tau})dt +
\sqrt{|V_t|}dW_{1,t},\\
d Y_t &= ({a_{21}^\lambda}V_t+a_{22} Y_t
+b_2-\lambda_2V_{t\wedge\tau})dt + \sqrt{|{V_t}|}d{W}_{2,t},
\end{align*}
under $\PP$ with time restricted to $[0,T]$. Therefore,
$((V,Y),W)$ is a solution to~(\ref{eq:v}),~(\ref{eq:y}) under
$\PP$, when time is restricted to $[0,\tau\wedge T]$. By
equivalence of $\PP$ and $\Q^\lambda$, we have $\PP(\tau<t)>0$ for
all {$t\in(0,T]$}.
\end{proof}}

\appendix
\section{Measure transformation for $\mathbb{A}_m(p)$}\label{sec:novfel}
\noindent In this section we prove that the exponential process
$L^\lambda_t$ defined by~(\ref{eq:den}) is a martingale for all
$\lambda\in\R^p$ for the class $\mathbb{A}_m(p)$. We present two
methods. The first is by verifying Novikov's condition by making
use of the explicit form of the square root SDE. The second method
uses pathwise uniqueness and also applies to a more general
situation.
\subsection{Using Novikov's condition}
As mentioned in~\cite{cher2007} page~129, a local version of
Novikov's condition holds for square root SDEs which satisfy the
Feller conditions. A good reference is lacking though. The
1-dimensional case, equivalent to the Cox-Ingersoll-Ross model, is
treated in~\cite{wong2006} and the proof uses an application of
the Feynman-Kac formula. For the general case $\mathbb{A}_m(p)$ we
will present a different method to verify Novikov's condition. The
underlying idea has also been used for verifying Novikov's
condition for the class $\mathbb{C}(2)$ without the Feller
conditions, see Section~\ref{sec:nov}. Note that under the Feller
conditions $V_{i,t}\geq 0$ almost surely, so $L^\lambda_t$ defined
by~(\ref{eq:den}) can be written as
\begin{equation}\label{eq:den3}
L^\lambda_t=\mathcal{E}(\int_0^\cdot\lambda^\top\sqrt{v(X_s)}d
W_s)_t.
\end{equation}
We prove that a local version of Novikov's condition holds for
$\mathbb{A}_m(p)$ under the additional requirement
\begin{equation}\label{eq:addreq}
\mbox{$\exists c_i>0$ for $i\leq m$ such that
$\sum_{j}c_ja_{ji}\leq -\half mc_i^2$ for all $i\leq m$.}
\end{equation}
For $m=p=2$, elementary but tedious computations show
that~(\ref{eq:addreq}) is satisfied in the following four cases:
\[
\begin{array}{rlll}
{\rm (i)}   &  a_{12}, a_{21}\geq 0, & a_{11}, a_{22}<0,  & \det a >0 \\
{\rm (ii)}  & a_{12}, a_{21}<0, & a_{11}, a_{22}\geq 0,  & \det a <0 \\
{\rm (iii)} &  a_{11},a_{12}<0\\
{\rm (iv)}  &  a_{22},a_{21}<0.
\end{array}
\]
Notice that the first case involves a sharpening of the weak
Feller conditions for $\mathbb{C}_2(2)$ on the elements of $a$.
This illustrates that Condition~(\ref{eq:addreq}) is not vacuous
and not in contradiction with the Feller conditions.

\begin{prop}\label{prop:nov} Consider a {solution $(X,W)$
on some filtered probability space
$(\Omega,\mathcal{F},(\mathcal{F}_t),\PP)$ to a} square root SDE
from $\mathbb{A}_m(p)$ with parameters $a,b,\alpha,\beta$. Assume
in addition that~(\ref{eq:addreq}) holds. Then for all
$\lambda\in\R^p$ Novikov's condition is satisfied for
$L^\lambda_t$ defined by~(\ref{eq:den3}), whence $L^\lambda_t$ is
a martingale for all $\lambda\in\R^p$.
\end{prop}
\begin{proof}
Let $\lambda\in\R^p$ be arbitrary. It is sufficient to find
$\varepsilon>0$ such that for all $t\geq0$ we have $
\E\exp(\half\int_t^{t+\varepsilon}\lambda^\top v(X_s)\lambda
ds)<\infty.$ Note that by the canonical representation and by the
Feller conditions~(\ref{eq:canfel}), this expectation reduces to
the form
\[
\E\exp(\int_t^{t+\varepsilon}(q_0+\sum_{i=1}^m q_i X_{i,s})
ds),\mbox{ for some $q_j\geq0$, $j=0,1,\ldots,m$.}
\]
Take $\varepsilon>0$ such that $\varepsilon q_i\leq c_i$ for all
$i\leq m$. Since $X$ solves~(\ref{eq:smp}) in canonical form
(see~(\ref{eq:can})), one gets
\begin{align*}
&\E\exp(\int_t^{t+\varepsilon}\sum_{i=1}^m q_i X_{i,s}
ds)\leq\frac{1}{\varepsilon}\int_t^{t+\varepsilon}
\E\exp(\sum_{i=1}^m c_i X_{i,s}) ds\\
&=\frac{1}{\varepsilon}\int_t^{t+\varepsilon}\exp(\sum_{i=1}^m
c_i(x_{i,0}+b_i s) \\
&\qquad\qquad\E\exp(\sum_{i=1}^m (\int_0^s \sum_{j=1}^m c_i
a_{ij}X_{j,u}du+\int_0^s c_i\sqrt{X_{i,u}}dW_{i,u}))ds.
\end{align*}
Interchanging the summation indices, applying a general form of
H\"{o}lder's inequality and using the assumptions on $c_i$, we
obtain that
\begin{align*}
&\E\exp(\sum_{i=1}^m (\int_0^s \sum_{j=1}^m
c_ia_{ij}X_{j,u}du+\int_0^s
c_i\sqrt{X_{i,u}}dW_{i,u}))\\
&\leq\prod_{i=1}^m[\E\exp(\int_0^s \sum_{j=1}^m mc_j
a_{ji}X_{i,u}du+\int_0^s
mc_i\sqrt{X_{i,u}}dW_{i,u})]^{1/m}\\
&\leq\prod_{i=1}^m[\E\exp(\int_0^s -\half
(mc_i)^2X_{i,u}du+\int_0^s
mc_i\sqrt{X_{i,u}}dW_{i,u})]^{1/m}\\
&=\prod_{i=1}^m[\E\mathcal{E}(mc_i\sqrt{X_{i}}\cdot
W_{i})]^{1/m}\leq1,
\end{align*}
where the last inequality holds by the supermartingale property of
an exponential process. The result follows.
\end{proof}
\medskip
When the additional requirement~(\ref{eq:addreq}) does not hold,
we cannot verify Novikov's condition for proving that
$L_t^\lambda$ is a martingale. However, applying the above
proposition twice solves this problem. We first transform the SDE
such that~(\ref{eq:addreq}) does hold and then transform it back
to the desired SDE.  This is possible since the above proposition
is valid for all $\lambda\in\R^p$.
\begin{prop} Consider a {solution $(X,W)$
on some filtered probability space
$(\Omega,\mathcal{F},(\mathcal{F}_t),\PP)$ to a} square root SDE
from $\mathbb{A}_m(p)$ with parameters $a,b,\alpha,\beta$. Then
$L^\lambda_t$ defined by~(\ref{eq:den3}) is a martingale for all
$\lambda\in\R^p$.
\end{prop}
\begin{proof} Let $c_i>0$ be arbitrary, $i\leq m$.
It is possible to choose $\mu\in\R^p$ such that
\[
c_i(a_{ii}+\mu_i)+\sum_{j\not=i}c_ja_{ji}\leq -\half mc_i^2,\mbox{
for }i\leq m,\mbox{ and }\mu_i=0\mbox{ for }i>m.
\]
We first show that $L^\mu_t=\mathcal{E}(\mu^\top\sqrt{v(X)}\cdot
W)_t$ is a martingale. Therefore, we consider the SDE in
$\mathbb{A}_m(p)$ with parameters $a^\mu,b,\alpha,\beta$, with
$a^\mu_{ii}=a_{ii}+\mu_i$, for $i\leq m$, $a^\mu_{ij}=a_{ij}$
otherwise. This SDE satisfies the conditions of
Proposition~\ref{prop:nov}. By existence of a (strong or weak)
solution, there exists a {filtered} probability space
$(\widehat{\Omega},\widehat{\mathcal{F}},{(\widehat{\mathcal{F}}_t)},\widehat{\Q}^\mu)$
with an {adapted} process $\widehat{X}$ and a
$\widehat{\Q}^\mu$-Brownian motion $\widehat{W}^\mu$ such that
\[
d \widehat{X}_t = (a^\mu\widehat{X}_t+b)dt
+\sqrt{v(\widehat{X}_t)}d\widehat{W}^\mu_t.
\]
By Proposition~\ref{prop:nov} the exponential process
\[
\widehat{L}^{-\mu}_t:=\mathcal{E}(-\mu^\top\sqrt{v(\widehat{X})}\cdot
\widehat{W}^\mu)_t
\]
is a $\widehat{\Q}^\mu$-martingale. Moreover, for a fixed
arbitrary $T>0$, the stopped process $\widehat{L}^{-\mu}_{T\wedge
t}$ is uniformly integrable. Hence we can change the measure
$\widehat{\Q}^\mu$ into an equivalent measure $\widehat{\PP}$ {on
$\mathcal{F}_\infty$} by
{$d\widehat{\PP}=\widehat{L}^{-\mu}_Td\widehat{\Q}^\mu$}. Then
$\widehat{W}$ defined by
\[
d\widehat{W}_t=d \widehat{W}^\mu_t+\mu^\top
\sqrt{v(\widehat{X}_t)}1_{\{t\leq T\} }dt,
\]
is a $\widehat{\PP}$-Brownian motion. Furthermore,
$(\widehat{X},\widehat{W})$ is a second solution to the initial
square root SDE
\[
d \widehat{X}_t = (a \widehat{X}_t+b)dt
+\sqrt{v(\widehat{X}_t)}d\widehat{W}_t,\mbox{
$\widehat{\PP}$-a.s.,}
\]
with time restricted to $[0,T]$. Hence applying
Proposition~\ref{prop:stop} (with $A=\emptyset$) gives that $X^T$
and $\widehat{X}^T$ have the same law. Moreover, for $t\leq T$ we
have
\[
(\widehat{L}^{-\mu}_t)^{-1}=\mathcal{E}(\mu^\top\sqrt{v(\widehat{X})}\cdot
\widehat{W})_t=\mathcal{E}(\mu^\top(\widehat{X}-\widehat{X}_0-\int_0^\cdot(a\widehat{X}_s+b)ds))_t,\,\widehat{\PP}\mbox{-a.s.,}
\]
and
\[
{L}^{\mu}_t=\mathcal{E}(\mu^\top\sqrt{v({X})}\cdot
{W})_t=\mathcal{E}(\mu^\top(X-X_0-\int_0^\cdot(aX_s+b)ds))_t,\,{\PP}\mbox{-a.s.,}
\]
so $(\widehat{L}^{-\mu}_t)^{-1}$ and ${L}^{\mu}_t$ are equal in
law for $t\leq T$. By equivalence of $\widehat{\PP}$ and
$\widehat{\Q}^\mu$, it holds that $\E_{\widehat{\PP}}
(\widehat{L}^{-\mu}_t)^{-1}=1$ for all $t\geq 0$. Therefore, $\EP
L^\mu_t=1$ for $t\leq T$. In fact, $\EP L^\mu_t=1$ for all $t\geq
0$, as $T$ can be chosen arbitrarily, whence $L^\mu_t$ is a
$\PP$-martingale.

Now we show that $L^\lambda_t$ is a $\PP$-martingale. Again fix an
arbitrary $T>0$. We change the measure ${\PP}$ into an equivalent
measure $\Q^\mu$ {on $\mathcal{F}_T$} by $d\Q^\mu=L^\mu_Td\PP$ and
see that
\[
d {X}_t = (a^\mu{X}_t+b)dt +\sqrt{v(X_t)}dW^\mu_t,
\]
for $t\leq T$, with $W^\mu$ a $\Q^\mu$-Brownian motion on $[0,T]$
given by $dW^\mu_t=dW_t-\mu^\top\sqrt{v(X_t)}dt$. Applying
Proposition~\ref{prop:nov} again gives that
$\mathcal{E}(\nu^\top\sqrt{v(X)}\cdot W^\mu)$ is a
$\Q^\mu$-martingale for all $\nu\in\R^p$ on $[0,T]$, whence on
$[0,\infty)$ since $T>0$ is arbitrary. Choosing $\nu=\lambda-\mu$
gives that
\begin{align*}
\EP L^\lambda_t=\E_{\Q^\mu} L^\lambda_t
(L^\mu_t)^{-1}&=\E_{\Q^\mu}
\mathcal{E}(\lambda^\top\sqrt{v(X)}\cdot
W)_t\mathcal{E}(-\mu^\top\sqrt{v(X)}\cdot W^\mu)_t\\
&=\E_{\Q^\mu} \mathcal{E}((\lambda-\mu)^\top\sqrt{v(X)}\cdot
W^\mu)_t=1,
\end{align*}
for all $t\geq0$, which completes the proof.
\end{proof}
\subsection{Using pathwise uniqueness}
For verifying Novikov's condition one needs the explicit form of
the underlying SDE. There are more general results in the
literature for proving that an exponential process is a martingale
without using the parameters of the SDE explicitly, for example
those contained in~\cite{Kallsen}, which treats the problem for
Dol\'eans exponentials of affine semimartingales under Feller
conditions. These  results cannot directly be used to show that
the process defined in~(\ref{eq:den3}) is a martingale. Other
results are
 the theorems in~\cite{wong2004}, Theorem~1 in~\cite{cher2007} and Theorem~A.1 in~\cite{hest2007}. The latter
theorem generalizes Theorem~1 in~\cite{cher2007}, but is not
applicable for the square root SDEs in $\mathbb{A}_m(p)$, since
strictly positiveness of the diffusion part $\sigma$ is required.
Moreover, it only treats the one-dimensional case.  Therefore, we
give another generalization of Theorem~1 in~\cite{cher2007} which
also applies to $\mathbb{A}_m(p)$. The proof goes along the same
{line of thought}, but for clarity we give it again, also to
emphasize the need for existence and uniqueness of strong
solutions, an aspect not mentioned in~\cite{cher2007}. {Note that
the latter implies that we cannot apply the result to the general
class $\mathbb{S}_m(p)$, as it is not clear whether pathwise
uniqueness holds for general square root SDEs.}

{In the proof of the next proposition we use a uniqueness in law
result for two weak solutions to an SDE up to a stopping time.
This result is stated and proved in Proposition~\ref{prop:stop} in
the next section. Its proof uses the existence of strong
solutions.}
\begin{prop}\label{prop:change2} Let
$\mu:\R^p\rightarrow\R^p$, $\sigma:\R^p\rightarrow\R^{p\times p}$,
$\gamma:\R^p\rightarrow\R^p$. Suppose we have weak existence for
the $p$-dimensional SDE
\begin{equation}\label{eq:SDE1}
dX_t=\mu(X_t)dt+\sigma(X_t)dW_t,
\end{equation}
as well as {existence and uniqueness of a strong solution} for the
SDE
\begin{equation}\label{eq:SDE2}
dX_t=(\mu(X_t)+\sigma(X_t)\gamma(X_t))dt+\sigma(X_t)dW_t.
\end{equation}
{Let $(X,W)$ be a weak solution to~(\ref{eq:SDE1}) on some
filtered probability space
$(\Omega,\mathcal{F},(\mathcal{F}_t),\PP)$ and
$(\widehat{X},\widehat{W})$ a solution to~(\ref{eq:SDE2}) on a
(possibly different) probability space
$(\widehat{\Omega},\widehat{\mathcal{F}},(\widehat{\mathcal{F}}_t),\widehat{\PP})$.
Suppose $\gamma(X_t)$ and $\gamma(\widehat{X}_t)$ have continuous
sample paths under $\PP$ respectively $\widehat{\PP}$. Then
$Y_t=\mathcal{E}(\gamma(X_t)\cdot W)_t$ is a $\PP$-martingale.}
\end{prop}
\begin{proof} Fix $T>0$ arbitrarily and define for
$n\in\mathbb{N}$
\[
\tau_n=\inf\{t>0:\|\gamma(X_t)\|\geq n\}\wedge
T,\quad\widehat{\tau}_n=\inf\{t>0:\|\gamma(\widehat{X}_t)\|\geq
n\}\wedge T.
\]
Then $Y^n_t=\mathcal{E}(\gamma(X)1_{[0,\tau_n]}\cdot W)_t$ is a
$\PP$-martingale, since Novikov's condition holds. Furthermore,
$Y^n_{t\wedge T}=Y^n_t$, so $Y^n$ is uniformly integrable, whence
$Y^n_\infty$ exists and equals $Y^n_T$. We can change the measure
$\PP$ into an equivalent measure $\Q^n$ on $\mathcal{F}_\infty$ by
{$d\Q^n=Y^n_\infty d\PP$}. Then $W^n$ defined by
$dW^n_t=dW_t-Y^n_tdt$ is a $\Q^n$-Brownian motion and $(X,W^n)$
satisfies
\[
dX_t=(\mu(X_t)+\sigma(X_t)\gamma(X_t)1_{\{t\leq\tau_n\}})dt+\sigma(X_t)dW^n_t,
\]
under $\Q^n$. Therefore, $(X,W^n)$ is a solution
to~(\ref{eq:SDE2}) on
$(\Omega,\mathcal{F}_\infty,(\mathcal{F}_t),\Q^n)$ with time
restricted to $[0,\tau_n]$.

By continuity of $\gamma(X_t)$ under $\PP$, we have
$\PP(\|\gamma(X_{\tau_n})\|\geq n\mbox{ or } \tau_n=T)=1$, whence
$\Q^n(\|\gamma(X_{\tau_n})\|\geq n\mbox{ or } \tau_n=T)=1$, for
all $n$, by equivalence of $\PP$ and $\Q^n$. So we can apply
Proposition~\ref{prop:stop} and obtain that $\tau_n$ and
$\widehat{\tau}_n$ have the same distribution under $\Q^n$
respectively $\widehat{\PP}$.

By continuity of $\gamma(X_t)$ under $\PP$, we have
$\tau_n\uparrow T$, $\PP$-a.s., which implies that
$Y^n_t1_{\{t\leq\tau_n\}}=Y_t1_{\{t\leq\tau_n\}}\uparrow
Y_t1_{\{t\leq T\}}$, $\PP$-a.s. Hence we can apply the Monotone
Convergence Theorem and obtain for $t\leq T$ that
\[
\EP Y_t = \lim_{n\rightarrow\infty} \EP
Y^n_t1_{\{t\leq\tau_n\}}=\lim_{n\rightarrow\infty} \Q^n(t<\tau_n)=
\lim_{n\rightarrow\infty} \widehat{\PP}(t<\widehat{\tau}_n)=1.
\]
where the last equality holds since
$\widehat{\PP}(\widehat{\tau}_n\uparrow T)=1$, by continuity of
$\gamma(\widehat{X}_t)$ under $\widehat{\PP}$. Because $T>0$ was
chosen arbitrarily, $\EP Y_t =1$ holds for all $t\geq 0$, whence
$Y_t$ is a $\PP$-martingale.
\end{proof}
\begin{remark}
Note that the above proposition implies existence and uniqueness
of a strong solution for~(\ref{eq:SDE1}). Indeed, suppose
$(X^1,W)$ and $(X^2,W)$ are solutions to~(\ref{eq:SDE1}) on some
filtered probability space
$(\Omega,\mathcal{F}_\infty,(\mathcal{F}_t),\PP)$. Fix $T>0$. Then
$\widehat{\PP}$ defined by $d\widehat{\PP}=Y_T d\PP$ is a
probability measure on $\mathcal{F}_T$, equivalent to $\PP$.
Furthermore, $\widehat{W}$ defined by $d\widehat{W}_t=dW_t-Y_tdt$
is a $\widehat{\PP}$-Brownian motion on $[0,T]$ and
$(X^1,\widehat{W})$ and $(X^2,\widehat{W})$ are solutions
to~(\ref{eq:SDE2}) up to $T$. By pathwise uniqueness
for~(\ref{eq:SDE2}) it holds that
\[
\widehat{\PP}(X^1_t=X^2_t,t\in[0,T])=1,
\]
which implies $\PP(X^1_t=X^2_t,t\in[0,T])=1$, as $\PP$ and
$\widehat{\PP}$ are equivalent. Since $T$ was chosen arbitrarily,
it follows that $\PP(X^1_t=X^2_t,\forall t\geq0)=1$. Hence
pathwise uniqueness holds for~(\ref{eq:SDE1}), which implies
existence and uniqueness of a strong solution, by Theorem~IV.2.1
in~\cite{Ikwat}.
\end{remark}
\begin{cor}
Consider a {solution $(X,W)$ on some filtered probability space
$(\Omega,\mathcal{F},(\mathcal{F}_t),\PP)$ to a} square root SDE
from $\mathbb{A}_m(p)$ with parameters $a,b,\alpha,\beta$. Then
$L^\lambda_t$ defined by~(\ref{eq:den3}) is a martingale for all
$\lambda\in\R^p$.
\end{cor}
\begin{proof}
This follows from Proposition~\ref{prop:change2} with
$\gamma(X_t)=\sqrt{v(X_t)}\lambda$. Both SDEs from the proposition
belong to $\mathbb{A}_m(p)$, for which weak existence holds by
continuity of the parameters and satisfaction of a {growth
condition}. Pathwise uniqueness holds by Theorem~1
in~\cite{YamadaI}, {which implies existence and uniqueness of a
strong solution, by Theorem~IV.2.1 in~\cite{Ikwat}.}
\end{proof}
\section{Uniqueness for a stopped SDE}
\noindent{As mentioned in the remark preceding
Proposition~\ref{prop:change2}, in this section we state and prove
a uniqueness in law result for two weak solutions to an SDE,
possibly defined on different probability spaces, up to a stopping
time. This result is stated in Proposition~\ref{prop:stop}. In the
proof we use a measurability lemma, which we prove first in
Lemma~\ref{lem:measstop}.}
\begin{lemma}\label{lem:measstop}
Let $(\Omega,\mathcal{F},(\mathcal{F}_t),\PP)$ be a filtered
probability space, $\tau$ a finite stopping time and $X$ a
stochastic process with continuous sample paths. If
$X_t:\Omega\rightarrow\R$ is $\mathcal{F}_{\tau+t}$-measurable for
all $t$, then $X_{t-\tau}1_{\tau< t}$ is
$\mathcal{F}_{t}$-measurable for all $t$.
\end{lemma}
\begin{proof}
It is possible to choose a sequence of stopping times
$\tau_n\downarrow\tau$ a.s.\ such that $\tau_n$ only assumes
countably many values. Since $ X_{t-\tau_n}1_{\tau_n< t}$
converges to $X_{t-\tau}1_{\tau< t}$ a.s., it is enough to prove
the statement for $\tau_n$ instead of $\tau$. For arbitrary Borel
set $B$ it holds that
\begin{align*}
\{X_{t-\tau_n}1_{\tau_n< t}\in B\}&=\{X_{t-\tau_n}\in B,\tau_n<
t\}\cup\{0\in B,\tau_n\geq
t\}\\
&=\bigcup_{k< t}\{X_{t-k}\in B,\tau_n=k\}\cup\{0\in B,\tau_n\geq
t\}\\
&=\bigcup_{k<t}(\{X_{t-k}\in B,\tau\leq
k\}\cap\{\tau_n=k\})\cup\{0\in B,\tau_n\geq t\}.
\end{align*}
Since $\tau_n$ is a stopping time, we have $\{0\in B,\tau_n\geq
t\}\in\mathcal{F}_t$ as well as $\{\tau_n=k\}\in\mathcal{F}_t$ for
$k<t$. Moreover, $X_t$ is $\mathcal{F}_{\tau+t}$-measurable for
all $t$, which means that
\[
\{X_t\in B\}\cap\{\tau+t\leq s\}\in\mathcal{F}_s, \quad\mbox{ for
all $t$ and $s$}.
\]
Choosing $s=t$ and substituting $t-k$ for $t$ in the above display
gives
\[
\{X_{t-k}\in B\}\cap\{\tau\leq k\}\in\mathcal{F}_t,
\]
which completes the proof.
\end{proof}
\begin{prop}\label{prop:stop}
Consider an SDE
\begin{equation}\label{eq:SDE}
dX_t = \mu(X_t)dt+\sigma(X_t)dW_t,
\end{equation}
which has a {unique strong solution $(X,W)$} on
$(\Omega,\mathcal{F},(\mathcal{F}_t),\PP)$. Let $\tau$ be a
stopping time of the form
\[
\tau = \inf\{t>0:X_t\in A\}\wedge {T},
\]
with $A$ a measurable set and {$T>0$}. {Let
$(\widehat{X},\widehat{W})$ be an adapted stochastic process on a
filtered probability space
$(\widehat{\Omega},\widehat{\mathcal{F}},(\widehat{\mathcal{F}}_t),\widehat{\PP})$,
with $\widehat{W}$ a $\widehat{\PP}$-Brownian motion. Suppose
$(\widehat{X},\widehat{W})$ is also a solution to~(\ref{eq:SDE})
under $\widehat{\PP}$, but on the stopped interval
$[0,\widehat{\tau}]$,} where we write
\[
\widehat{\tau} = \inf\{t>0:\widehat{X}_t\in A\}\wedge T.
\]
If $\widehat{\PP}(\widehat{X}_{\widehat{\tau}}\in A\mbox{ or }
\widehat{\tau}=T)=1$, then the stopping times $\tau$ and
$\widehat{\tau}$ as well as the stopped processes $X^\tau$ and
$\widehat{X}^{\widehat{\tau}}$ have the same distribution under
$\PP$ respectively $\widehat{\PP}$.
\end{prop}
\begin{proof}
We extend the solution $\widehat{X}$ to~(\ref{eq:SDE}) on
$[0,\widehat{\tau}]$ to a solution $Y$ to~(\ref{eq:SDE}) on the
whole interval $[0,\infty)$, for which we use existence of a
strong solution. Define a filtration $(\mathcal{G}_t)$ by
$\mathcal{G}_t:=\widehat{\mathcal{F}}_{\widehat{\tau} +t}$. Then
$\widetilde{W}_t:=\widehat{W}_{\widehat{\tau}+t}-\widehat{W}_{\widehat{\tau}}$
is a $\widehat{\PP}$-Brownian motion with respect to
$(\mathcal{G}_t)$. By existence of a strong solution, there exists
a process $Z$ adapted to $\mathcal{G}$ with initial value
$\widehat{X}_{\widehat{\tau}}$ such that
\[
Z_t=\widehat{X}_{\widehat{\tau}}+\int_0^t \mu(Z_s)ds+\int_0^t
\sigma(Z_s)d\widetilde{W}_s.
\]
Define
\[
Y_t=\widehat{X}_{t}1_{t\leq
\widehat{\tau}}+Z_{t-\widehat{\tau}}1_{t>\widehat{\tau}}.
\]
By Lemma~\ref{lem:measstop}, $Y_t$ is
$\widehat{\mathcal{F}}_t$-measurable. It holds that
$Y_{\widehat{\tau}+t}=Z_t$ for $t\geq0$, so that
\begin{align*}
Y_{\widehat{\tau}+t}=Z_t&=\widehat{X}_{\widehat{\tau}}+\int_0^t
\mu(Y_{\widehat{\tau}+s})ds+\int_0^t
\sigma(Y_{\widehat{\tau}+s})d(\widehat{W}_{\widehat{\tau}+s}-\widehat{W}_{\widehat{\tau}})\\
&=\widehat{X}_{\widehat{\tau}}+\int_{\widehat{\tau}}^{\widehat{\tau}+t}
\mu(Y_{s})ds+\int_{\widehat{\tau}}^{\widehat{\tau}+t}
\sigma(Y_{s})d\widehat{W}_{s}.
\end{align*}
Note that
\begin{align*}
\widehat{X}_{\widehat{\tau}}&=\widehat{X}_0+\int_0^{\widehat{\tau}}\mu(\widehat{X}_s)ds+\int_0^{\widehat{\tau}}\sigma(\widehat{X}_s)d\widehat{W}_s\\
&=\widehat{X}_0+\int_0^{\widehat{\tau}}\mu(Y_s)ds+\int_0^{\widehat{\tau}}\sigma(Y_s)d\widehat{W}_s,
\end{align*}
whence
\begin{align*}
Y_{t}1_{t>
\widehat{\tau}}=Y_{\widehat{\tau}+t-\widehat{\tau}}1_{t>\widehat{\tau}}&=(\widehat{X}_{\widehat{\tau}}+\int_{\widehat{\tau}}^{t}
\mu(Y_{s})ds+\int_{\widehat{\tau}}^{t}
\sigma(Y_{s})d\widehat{W}_{s})1_{t>\widehat{\tau}}\\
&=(Y_0+\int_{0}^{t} \mu(Y_{s})ds+\int_0^{t}
\sigma(Y_{s})d\widehat{W}_{s})1_{t>\widehat{\tau}}.
\end{align*}
On the other hand it holds that
\begin{align*}
Y_{t}1_{t\leq\widehat{\tau}}=\widehat{X}_{t}1_{t\leq
\widehat{\tau}}&=(\widehat{X}_0+\int_0^t\mu(\widehat{X}_s)ds+\int_0^t\sigma(\widehat{X}_s)d\widehat{W}_s)1_{t\leq\widehat{\tau}}\\
&=(Y_0+\int_0^t\mu(Y_s)ds+\int_0^t\sigma(Y_s)d\widehat{W}_s)1_{t\leq\widehat{\tau}}.
\end{align*}
Hence $(Y,\widehat{W})$ is a solution to~(\ref{eq:SDE}) on
$(\widehat{\Omega},\widehat{\mathcal{F}},(\widehat{\mathcal{F}}_t),\widehat{\PP})$.
By uniqueness in distribution, $Y$ and $X$ have the same
distribution. Since the paths of $\widehat{X}$ and $Y$ coincide
for $t\leq \widehat{\tau}$ and since
$\widehat{\PP}(\widehat{X}_{\widehat{\tau}}\in A\mbox{ or }
\widehat{\tau}=T)=1$, it holds $\widehat{\PP}$-almost surely that
\begin{align*}
\widehat{\tau} &=\inf\{t>0:Y_t\in A\}\wedge T.
\end{align*}
Comparing this with the expression for $\tau$, we see that
$\widehat{\tau}$ and $\tau$ as well as $Y^{\widehat{\tau}}$ and
$X^{\tau}$ have the same distribution. From
$Y^{\widehat{\tau}}=\widehat{X}^{\widehat{\tau}}$ it follows that
$\widehat{X}^{\widehat{\tau}}$ and $X^{\tau}$ have the same
distribution.
\end{proof}
\bibliographystyle{plain}
\bibliography{refs}

\end{document}